\documentclass{article}
\usepackage[utf8]{inputenc}

\usepackage[english]{babel}

\usepackage[utf8]{inputenc}
\usepackage[T1]{fontenc}

\usepackage{amssymb}
\usepackage{stmaryrd}
\usepackage{amsmath,amsfonts}
\usepackage[tikz]{bclogo}
\usepackage[top=4.34cm]{geometry}

\usepackage{tikz}
\usepackage{tkz-tab}
\usepackage{caption}
\usepackage{latexsym}
\usepackage{amssymb}
\usepackage{amsmath}
\usepackage{subfig}
\usetikzlibrary{shapes.geometric}
\usetikzlibrary{decorations.pathreplacing}

\usepackage{lmodern}
\usepackage{textcomp}
\usepackage{mathtools, bm}
\usepackage{amssymb, bm}
\usepackage[unq]{unique}
\usepackage{enumerate}
\usepackage{comment}
\usepackage[breaklinks]{hyperref}

\usepackage[numbers]{natbib}  

\usepackage[breaklinks]{hyperref}
\usepackage[numbers]{natbib}
\usepackage{graphicx} 

\usepackage{amsthm}
\usepackage{multirow}

\usepackage{cleveref}

\newtheorem{theorem}{Theorem}[section]
\newtheorem{lemma}[theorem]{Lemma}

\newtheorem{conjecture}[theorem]{Conjecture}

\newtheorem{Claim}[theorem]{Claim}
\newtheorem{remark}[theorem]{Remark}
\newtheorem{Corollary}[theorem]{Corollary}
\newtheorem{strategy}[theorem]{Strategy}

\newcommand{\cM}{\mathcal{M}}
\newcommand{\cP}{\mathcal{P}}

\newcommand{\distR}{\text{dist}_R}

\DeclareMathOperator{\Tr}{T_1}
\DeclareMathOperator{\So}{S_0}
\DeclareMathOperator{\Sn}{S_1}
\DeclareMathOperator{\Co}{C_0}
\DeclareMathOperator{\Cn}{C_1}
\DeclareMathOperator{\Oo}{O_0}
\DeclareMathOperator{\On}{O_1}

\DeclareMathOperator{\Dn}{D_1}
\DeclareMathOperator{\Dt}{D_{2^+}}

\title{The asymptotic of off-diagonal online Ramsey numbers for paths}
\author{Adva Mond \footnote{\href{mailto:am2759@cam.ac.uk}{am2759@cam.ac.uk}, Department of Pure Mathematics and Mathematical Statistics (DPMMS), University of Cambridge, Wilberforce Road, Cambridge, CB3 0WA, United Kingdom} \and Julien Portier \footnote{\href{mailto:jp899@cam.ac.uk}{jp899@cam.ac.uk}, Department of Pure Mathematics and Mathematical Statistics (DPMMS), University of Cambridge, Wilberforce Road, Cambridge, CB3 0WA, United Kingdom}}
\date{}

\begin{document}

\maketitle

\begin{abstract}
    We prove that for every $k\ge 10$, the online Ramsey number for paths $P_k$ and $P_n$ satisfies $\tilde{r}(P_k,P_n) \geq \frac{5}{3}n + \frac{k}{9} - 4$, matching up to a linear term in $k$ the upper bound recently obtained by Bednarska-Bzd{\k{e}}ga [European Journal of Combinatorics 118 (2024)].
    In particular, this implies $\lim_{n \rightarrow \infty} \frac{\tilde{r}(P_k, P_n)}{n} = \frac{5}{3}$, whenever $10 \le k=o(n)$, disproving a conjecture by Cyman, Dzido, Lapinskas and Lo [Electronic Journal of Combinatorics (2015)].
\end{abstract}

\section{Introduction}

Online Ramsey numbers, sometimes referred to as online size Ramsey numbers, were firstly introduced by Beck \cite{beck1993achievement} and later rediscovered by Kurek and Ruci{\'n}ski \cite{kurek2005two}.
One way to define them is by a one-player game where the board is the infinite graph $K_{\mathbb N}$, with a $2$-colouring of its edges by red and blue, which is hidden from the player.
Given two (finite) graphs $G$ and $H$, at each round the player chooses an edge of the graph and its colour is revealed, and the game ends when either a red copy of $G$ or a blue copy of $H$ is revealed.
The online Ramsey number $\tilde{r}(G,H)$ is the minimum number of rounds needed for the player to end the game for any colouring of the board.
A more common but equivalent definition of the online Ramsey number $\tilde{r}(G,H)$ is by a combinatorial game $\tilde{R}(G,H)$ played by two Players, Builder and Painter, where the board is the infinite vertex set $\mathbb N$.
At each round, Builder adds a previously not yet selected edge between two vertices and Painter colours it red or blue.
The game ends when either a red copy of $G$ or a blue copy of $H$ is created.
Builder's goal is to end the game as fast as possible while Painter's goal is to make it last as long as possible.
The online Ramsey number $\tilde{r}(G,H)$ is the number of rounds in a game when both players play optimally.
This notion is closely related to the size Ramsey number $\hat{r}(G,H)$, the minimum number of edges in a graph for which every $2$-colouring of its edges by red and blue contains either a red copy of $G$ or a blue copy of $H$.
In particular, we have $\tilde{r}(G,H) \le \hat{r}(G,H)$.
If $G=H$ then it is common to write $\tilde{r}(G)$ and $\hat{r}(G)$ instead of $\tilde{r}(G,G)$ and $\hat{r}(G,G)$, respectively, for simplicity.

A significant attention was given to online Ramsey numbers for paths.
Beck~\cite{beck1983size} proved that the size Ramsey number $\hat{r}(P_n)$ is linear in $n$, implying that the online Ramsey number $\tilde{r}(P_n)$ is linear as well, where by $P_n$ we denote the path consisting of $n$ vertices.
Grytczuk, Kierstead and Pra\l at~\cite{grytczuk2008line}
and Pra\l at~\cite{pralat2008note,prlaat2012note} studied the value of $\tilde{r}(P_k,P_n)$ and determined it when $\max\{k,n \}\le 9$.
In~\cite{grytczuk2008line} they gave bounds also for the general case, showing that $k+n-3 \le \tilde{r}(P_k,P_n) \le 2k+2n-7$.
Recently, Bednarska-Bzd{\k{e}}ga~\cite{BEDNARSKABZDEGA2024103873} proved that $\tilde{r}(P_k,P_n) \leq \frac{5}{3}n+12k$, which improves the upper bound when $n$ is large compared to $k$.
As for the lower bound, Cyman, Dzido, Lapinskas and Lo~\cite{cyman2015line} proved that $\tilde{r}(P_k,P_n) \geq \frac{3}{2}n+\frac{k}{2}-\frac{7}{2}$ for every $k \geq 5$.
They believed that the strategy they suggested for Painter in their proof was asymptotically optimal, and therefore made the following conjecture.

\begin{conjecture}[\cite{cyman2015line}]
\label{conj:Czyman}
For every fixed $k \geq 5$, we have
    $$\lim_{n \rightarrow{} \infty} \frac{\tilde{r}(P_k,P_n)}{n}=\frac{3}{2}.$$
\end{conjecture}

Other variants of online Ramsey numbers for paths were also studied.
For example, in~\cite{balogh2022strengthening} and in~\cite{perez2021line} the authors studied online Ramsey numbers for ordered paths in infinite complete ordered graphs and hypergraphs, and in~\cite{blavzej2019induced} the authors considered the induced version.

In this paper we continue this line of research and study the online Ramsey number for paths $\tilde{r}(P_k,P_n)$.
The main result of our paper is a lower bound for $\tilde{r}(P_k,P_n)$ when $k \geq 10$, which matches the upper bound obtained by Bednarska-Bzd{\k{e}}ga 
up to a linear term in $k$, and therefore also disproves \Cref{conj:Czyman}.

\begin{theorem}
\label{thm:main}
    For any integers $n$ and $k\ge 10$ we have $$\tilde{r}(P_k,P_n) \geq \frac{5}{3} n + \frac{k}{9} - 4.$$
\end{theorem}
Combining this with the upper bound by Bednarska-Bzd{\k{e}}ga~\cite{BEDNARSKABZDEGA2024103873} mentioned above, we get the following result.
\begin{Corollary}
\label{cor:main}
    Let $10 \le k = o(n)$.
    Then we have
    $$\lim_{n \rightarrow{} \infty} \frac{\tilde{r}(P_k,P_n)}{n}=\frac{5}{3}.$$
\end{Corollary}

\subsection{Notation}
As mentioned above, for an integer $\ell \ge 1$ we denote by $P_{\ell}$ the path on $\ell$ vertices (and $\ell-1$ edges).
By a \emph{round} of the game we mean one turn of Builder followed by one turn of Painter.
The first round of the game is round $1$.
We often consider the graph spanned by the blue (red) edges which we call the \emph{blue (red) graph}.
We also say \emph{blue (red) component} to refer to a connected component in the blue (red) graph.

\section{Proof idea}
\label{sec:idea}

\Cref{thm:main} is in fact a corollary of the following two statements.
The first one gives a lower bound for the base case where $k=10$.
\begin{theorem}
\label{thm:k=10}
    For any integer $n$ we have
    \[\tilde{r}(P_{10}, P_n) \ge \frac{5}{3}n - 2. \]
\end{theorem}
\begin{lemma}
\label{lem:PkPm}
    Let $k \ge m$ be integers and $H$ be a graph.
    Then we have
    \[\tilde{r}(P_k,H) \ge \tilde{r}(P_m,H) + \left\lfloor\frac{k-1}{m-1} \right\rfloor - 1. \]
\end{lemma}
Taking $H = P_n$ and $m=10$, \Cref{lem:PkPm} together with \Cref{thm:k=10} immediately implies \Cref{thm:main}.
\begin{proof}[Proof of \Cref{thm:main} assuming \Cref{thm:k=10} and \Cref{lem:PkPm}]
Let $k\ge 10$.
Then
\begin{align*}
    \tilde{r}(P_k,P_n) &\ge \tilde{r}(P_{10},P_n) + \left\lfloor\frac{k-1}{9} \right\rfloor - 1 \ge \frac{5}{3}n + \frac{k}{9} - 4. \qedhere
\end{align*}
\end{proof}

The proofs of \Cref{thm:k=10} and \Cref{lem:PkPm} are given in \Cref{sec:potentialfunc}.
As the proof of \Cref{lem:PkPm} is quite straightforward, we give the main ideas for the proof of \Cref{thm:k=10}.

In order to prove a lower bound on $\tilde{r}(P_{10},P_n)$, we provide Painter with an explicit strategy (see \Cref{strategy:painter}).
Following the proofs of previous results, the first feature of our strategy is that it does not allow any red copy of $P_{10}$ to be created.
This is shown in \Cref{lem:NoRedPk}.
It is then left to prove that, if Painter follows this strategy, it takes at least $\frac{5}{3}n-2$ rounds for Builder to create a blue copy of $P_n$.

In our proof for \Cref{thm:k=10} we let Painter follow \Cref{strategy:painter}, and we denote the set of possible moves in this strategy by $\cM$.
We split all possible moves of Builder into categories, and for each category we instruct Painter with a response.
Additionally, for each move of Painter $M\in \cM$, we assign a variable $x_M$ to count the number of times this move was played throughout the game.
We then define a set $\cP$ of parameters of the board which we track throughout the course of the game.
These parameters are in fact functions of integers.
Given a parameter $X \in \cP$ we write $X(t)$ for the value of this parameter in the graph obtained after precisely $t$ rounds have been played, and accordingly we write $\cP(t) \coloneqq \left\{X(t) ~:~ X\in \cP \right\}$.
We denote by $N$ the number of rounds the game lasts until Builder wins.
In particular we can write $\sum_{M\in \cM} x_M = N$.

Two of the parameters we consider in $\cP$ play key roles in our proof.
One, denoted by $\Dn$, is the number of vertices of degree one in the blue graph, and the other, denoted by $\Dt$, is the number of vertices of degree at least two in the blue graph.
For example, if Painter coloured blue an edge added by Builder for the first time in the game in round $t_0$, then we have $\Dn(t_0) = 2$.
More importantly, since a red copy of $P_{10}$ is never created, at the end of the game there is a blue copy of $P_n$ and therefore we have $\Dn(N) + \Dt(N) \ge n$ and $\Dt(N) \ge n-2$.
Hence, our goal is to control the growth of the function $\Dn(t) + \Dt(t)$ throughout the game, showing that $\Dn(N) + \Dt(N) \ge n$ together with $\Dt(N) \ge n-2$ imply that $N$ is large enough.

To do so, we analyse the change of values of parameters in $\cP$ over the rounds, in terms of moves in $\cM$.
This yields a set of linear inequalities, where the final values of the parameters in $\cP(N)$ are bounded by linear combinations of the variables $\left\{x_M ~:~ M \in \cM \right\}$.
Recalling that $\sum_{M\in \cM} x_M = N$, and considering the set of linear inequalities mentioned above, we can bound $N$ from below by solving the linear programming problem associated.
In our specific case, this linear programming problem can be solved by a linear combination of the inequalities on $\cP(N)$, and we define this exact linear combination to be our \emph{potential function} $\beta$ (see (\ref{eq:beta})). Naturally, as $\beta$ is a linear combination of the parameters in $\cP$, it is also a function of $t$.
In fact, $\beta(t)$ is bounded from below by $\frac{2}{3}\Dn(t) + \frac{5}{3}\Dt(t)$.
Hence, we control the growth of $\Dn(t) + \Dt(t)$ by bounding the increment of the potential function $\beta$ over the rounds.

More precisely, our goal is to show that the potential function $\beta$ satisfies the following properties.
\begin{enumerate}
    \item[(1)] $\beta(0) = 0$.
    \item[(2)] $\beta(t) - \beta(t-1) \le 1$ for every $t \in [N]$.
    \item[(3)] $\beta(N) \ge \frac{5}{3}n - 2$.
\end{enumerate}
If $\beta$ does satisfy these properties then
\[N \ge \sum_{t\in [N]} \beta(t) - \beta(t-1) = \beta(N) - \beta(0) \ge \frac{5}{3}n - 2, \]
implying \Cref{thm:k=10}.
Property (1) will follow immediately from the definition of $\beta$.
Property (2) is where we bound the increment of $\beta$.
We prove this in \Cref{lem:IncrementPotential} which is our main lemma.
Then property (3) will follow easily from $\beta(t) \ge \frac{2}{3}\Dn(t) + \frac{5}{3}\Dt(t)$ as we have $\Dn(N) + \Dt(N) \ge n$ and $\Dt(N) \ge n-2$.
Hence, \Cref{lem:IncrementPotential} is the heart of our argument.

We would like to note that the use of the potential function $\beta$ is not completely needed in our proof, as one could solve the linear programming problem associated and derive the same bound.
However, this would lead to a tedious computation involving many parameters, whereas when using the potential function $\beta$ it suffices to bound its increments case by case.

The reason that \Cref{strategy:painter} gives the asymptotically correct lower bound lies in a certain core idea implemented in Painter's moves.
To point it out, we give some background on the proofs of the upper and lower bounds for $\tilde{r}(P_k,P_n)$ as given in~\cite{BEDNARSKABZDEGA2024103873,cyman2015line}.
In the proof of the upper bound in~\cite{BEDNARSKABZDEGA2024103873}, Builder follows a strategy which starts by constructing $\frac{1}{3}n + o(n)$ many blue copies of $P_3$, which lasts for $n+o(n)$ rounds of the game.
This is Stage $1$ of the strategy.
Then, in Stages $2$ and $3$, Builder merges all these blue paths into one blue copy of $P_n$.
For these two stages $\frac{2}{3}n + o(n)$ more rounds are enough.
This results in a total of $\frac{5}{3}n + o(n)$ rounds, and shows that $\tilde{r}(P_k,P_n) \le \frac{5}{3}n + o(n)$.

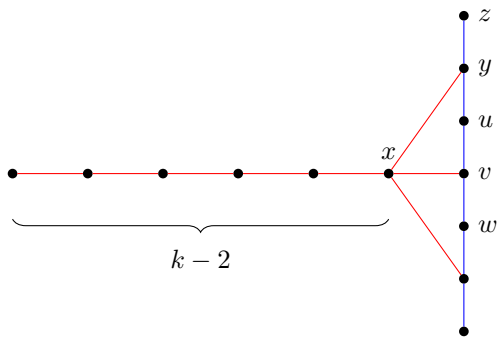
\begin{figure}
\captionsetup[subfigure]{labelformat=empty}
\centering
\begin{tikzpicture}[auto, vertex/.style={circle,draw=black!100,fill=black!100, thick,inner sep=0pt,minimum size=1mm},  novertex/.style={circle,draw=white!100}]
	\node (x1) at (-3,0) [vertex] {};
    \node (x2) at (-2,0) [vertex] {};
    \node (x3) at (-1,0) [vertex] {};
    \node (x4) at (0,0) [vertex] {};
    \node (x5) at (1,0) [vertex] {};
    \node (x6) at (2,0) [vertex,label=$x$] {};
    \node (x7) at (3,0) [vertex,label=right:$v$] {};
    \node (x8) at (3,1.4) [vertex,label=right:$y$] {};
    \node (x9) at (3,-1.4) [vertex] {};
    \node (y1) at (3,2.1) [vertex,label=right:$z$] {};
    \node (y2) at (3,0.7) [vertex,label=right:$u$] {};
    \node (y3) at (3,-0.7) [vertex,label=right:$w$] {};
    \node (y4) at (3,-2.1) [vertex] {};
		
	\draw [-,red] (x1) --node[inner sep=2pt,swap]{} (x2);
    \draw [-,red] (x2) --node[inner sep=2pt,swap]{} (x3);  
    \draw [-,red] (x3) --node[inner sep=2pt,swap]{} (x4);
    \draw [-,red] (x4) --node[inner sep=2pt,swap]{} (x5);
    \draw [-,red] (x5) --node[inner sep=2pt,swap]{} (x6);
    \draw [-,red] (x6) --node[inner sep=2pt,swap]{} (x7);
    \draw [-,red] (x6) --node[inner sep=2pt,swap]{} (x8);
    \draw [-,red] (x6) --node[inner sep=2pt,swap]{} (x9);

    \draw [decorate,decoration={brace,amplitude=5pt,mirror,raise=4ex}]
  (-3,0) -- (2,0) node[midway,yshift=-4em]{$k-2$};

	\draw [-,blue] (y1) --node[inner sep=2pt,swap]{} (x8);
    \draw [-,blue] (x8) --node[inner sep=2pt,swap]{} (y2);
    \draw [-,blue] (y2) --node[inner sep=2pt,swap]{} (x7);
    \draw [-,blue] (x7) --node[inner sep=2pt,swap]{} (y3);
    \draw [-,blue] (y3) --node[inner sep=2pt,swap]{} (x9);
    \draw [-,blue] (x9) --node[inner sep=2pt,swap]{} (y4);
	\end{tikzpicture}
\caption{An example of a construction obtained when Painter follows the suggested strategy in~\cite{cyman2015line} and Builder plays such that the game lasts for $\frac{3}{2}n + o(n)$ rounds.}
\label{fig:outline1}
\end{figure}

In~\cite{cyman2015line} the authors prove the lower bound $\tilde{r}(P_k, P_n) \ge \frac{3}{2}n + o(n)$, for $k\ge 5$, by providing Painter with the following strategy.
Painter colours red every edge added by Builder, unless it creates a red copy of either a $P_k$ or a cycle, in which case Painter colours it blue.
This way, no red copy of $P_k$ can ever appear, so Builder wins after creating a blue copy of $P_n$.
By a careful analysis of this process, they prove that it takes Builder at least $\frac{3}{2}n+o(n)$ rounds to win.
Moreover, their analysis of this strategy is asymptotically optimal, as Builder can indeed create a blue copy of $P_n$ in $\frac{3}{2}n+o(n)$ rounds when Painter follows their suggested strategy.
Builder firstly creates a red copy of $P_{k-1}$, and then adds two edges $vu$ and $vw$ adjacent to one of its endpoints $v$, which are both coloured blue by Painter, as illustrated in \Cref{fig:outline1}.
This constructs a blue path $uvw$ of length $3$.
From this point on, Builder increases the length of this blue path by $2$ every three rounds.
In the next round Builder adds the edge $xy$, where $x$ is the neighbour on the red copy of $P_{k-1}$ of its endpoint $v$, and Painter colours it red.
This creates another red copy of $P_{k-1}$ which intersects an existing one with $k-2$ vertices.
Then Builder repeats the same moves played for the first red copy of $P_{k-1}$, building two edges $yu$ and $yz$.
These two edges are coloured blue by Painter, increasing the length of the blue path by $2$ within three rounds.
Repeating this (as illustrated in \Cref{fig:outline1} with one more repetition), Builder indeed wins after $\frac{3}{2}n + o(n)$ rounds.

The main reason that Painter's strategy from~\cite{cyman2015line} does not yield a lower bound greater than $\frac{3}{2}n + o(n)$ is that red copies of $P_{k-1}$ are easily built, and once a red copy of $P_{k-1}$ is built, Painter has no choice but colouring blue any edge incident to an endpoint of it.
This makes the process of building a blue copy of $P_n$ rather quick; as explained above, Builder can increase the length of a blue path by $2$ every three rounds.
Thus, the first idea of our proof is to prevent Builder from playing this sequence of moves.
One possible way for Painter to achieve this, is by colouring blue any edge $e=xy$ where $x$ is a second to last vertex on a red $P_{k-1}$.
Indeed, if $y$ is not adjacent to any blue edge and Painter colours the edge $e=xy$ blue, then Builder does not increase lengths of blue paths faster than Stage $1$ in the proof of the upper bound in~\cite{BEDNARSKABZDEGA2024103873}.
However, if $y$ is already adjacent to a blue edge, colouring it blue increases the length of a blue path too fast.
Therefore, in this case, we call the vertex $y$ \emph{terminal} and we instruct Painter to colour the edge red.
If Builder later claims another edge incident to $y$, Painter will be forced to colour it blue.
This results, in the worst-case, with the merger of two blue components in two moves, which is similar to Stages $2$ and $3$ in the proof of the upper bound in~\cite{BEDNARSKABZDEGA2024103873}.
In total, whether Painter colours the edge $xy$ blue or red, it is guaranteed that Builder will not increase lengths of blue paths too fast.
Therefore, it makes sense to instruct Painter to colour the edge $xy$ blue whenever $x$ is a non-endpoint vertex on a red path, except for when $y$ is already adjacent to one blue edge, in which case Painter colours it red.

\begin{figure}
\captionsetup[subfigure]{labelformat=empty}
\centering
\begin{tikzpicture}[auto, vertex/.style={circle,draw=black!100,fill=black!100, thick,inner sep=0pt,minimum size=1mm},  novertex/.style={circle,draw=white!100}]
	\node (x1) at (-3,0) [vertex] {};
    \node (x2) at (-2,0) [vertex] {};
    \node (x3) at (-1,0) [vertex] {};
    \node (x4) at (0,0) [vertex] {};
    \node (x5) at (1,0) [vertex] {};
    \node (x6) at (2,0) [vertex] {};
    \node (x7) at (3,0) [vertex] {};
    \node (y1) at (-3,1) [vertex] {};
    \node (y2) at (-2,1) [vertex] {};
    \node (y3) at (-1,1) [vertex] {};
    \node (y4) at (0,1) [vertex] {};
    \node (y5) at (1,1) [vertex] {};
    \node (y6) at (2,1) [vertex] {};
    \node (y7) at (3,1) [vertex] {};
    \node (z1) at (-3,-1) [vertex] {};
    \node (z2) at (-2,-1) [vertex] {};
    \node (z3) at (-1,-1) [vertex] {};
    \node (z4) at (0,-1) [vertex] {};
    \node (z5) at (1,-1) [vertex] {};
    \node (z6) at (2,-1) [vertex] {};
    \node (z7) at (3,-1) [vertex] {};
		
	\draw [-,red] (x1) --node[inner sep=2pt,swap]{} (x2);
    \draw [-,red] (x2) --node[inner sep=2pt,swap]{} (x3);  
    \draw [-,red] (x3) --node[inner sep=2pt,swap]{} (x4);
    \draw [-,red] (x4) --node[inner sep=2pt,swap]{} (x5);
    \draw [-,red] (x5) --node[inner sep=2pt,swap]{} (x6);
    \draw [-,red] (x6) --node[inner sep=2pt,swap]{} (x7);

	\draw [-,blue] (x1) --node[inner sep=2pt,swap]{} (y1);
    \draw [-,blue] (x1) --node[inner sep=2pt,swap]{} (z1);
    \draw [-,blue] (x2) --node[inner sep=2pt,swap]{} (y2);
    \draw [-,blue] (x2) --node[inner sep=2pt,swap]{} (z2);
    \draw [-,blue] (x3) --node[inner sep=2pt,swap]{} (y3);
    \draw [-,blue] (x3) --node[inner sep=2pt,swap]{} (z3);
    \draw [-,blue] (x4) --node[inner sep=2pt,swap]{} (y4);
    \draw [-,blue] (x4) --node[inner sep=2pt,swap]{} (z4);
    \draw [-,blue] (x5) --node[inner sep=2pt,swap]{} (y5);
    \draw [-,blue] (x5) --node[inner sep=2pt,swap]{} (z5);
    \draw [-,blue] (x6) --node[inner sep=2pt,swap]{} (y6);
    \draw [-,blue] (x6) --node[inner sep=2pt,swap]{} (z6);
    \draw [-,blue] (x7) --node[inner sep=2pt,swap]{} (y7);
    \draw [-,blue] (x7) --node[inner sep=2pt,swap]{} (z7);
	\end{tikzpicture}
\caption{A red copy of $P_7$ with blue copies of $P_3$ on each vertex of it.}
\label{fig:outline2}
\end{figure}
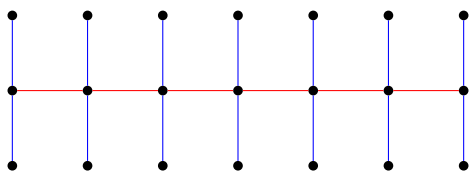

This idea alone is sufficient for disproving \Cref{conj:Czyman}.
Indeed, when $k$ is large enough, by a careful analysis of the strategy described above, it can be shown that there exists an absolute constant $\varepsilon > 0$ such that Builder cannot create a blue copy of $P_n$ in less than $\left( \frac{3}{2}+\varepsilon \right)n+o(n)$ rounds.
However, this is not sufficient for proving an asymptotically optimal lower bound of $\frac{5}{3}n + o(n)$ as in \Cref{thm:main}, since Builder can play as follows.
By building copies of $P_3$ centred in each vertex of a red copy of, say, $P_7$, Builder creates a structure as in \Cref{fig:outline2}.
Doing this, Builder creates $7$ blue copies of $P_3$ in $20$ rounds, which is slightly faster than Stage $1$ of the upper bound in~\cite{BEDNARSKABZDEGA2024103873}, where Builder creates an average of $1$ blue $P_3$ every $3$ rounds.
Following Stages $2$ and $3$ from the same proof, Builder can create a blue copy of $P_n$ within $\left(\frac{5}{3}-\delta \right)n+o(n)$ rounds, for $\delta=\frac{5}{63}$.
Thus, to prove the asymptotically optimal bound $\frac{5}{3}n+o(n)$ from \Cref{thm:main}, we must address this possible scenario in Painter's strategy.

Note that vertices which are not close to being ends of red copies of $P_{k-1}$ can have an incident red edge without causing the problem we mention above.
Hence, we wish to allow Painter to distinguish between those edges incident to vertices of red paths which can be coloured red and those that should be coloured blue.
To do so, in each large enough red component, three vertices which are not close to an end of any red copy of $P_{k-1}$ will be called \emph{central} vertices.
Having these central vertices, Painter can colour red the edges incident to them added by Builder, so that structures as the one in \Cref{fig:outline2} are prevented, while the main idea of the strategy is preserved.
For computational reasons, we need three central vertices per large red component, and not any less.

Our strategy for Painter does allow red copies of $P_9$ to be created, but not red copies of $P_{10}$.
This is the reason that our base case is $k=10$.

This summarises the ideas behind our proof of \Cref{thm:k=10}.
In \Cref{sec:strategy} we describe an explicit strategy for Painter and we state several immediate properties of it.
In \Cref{sec:potentialfunc} we define the set of parameters $\cP$ which we track throughout the course of the game, and we define the potential function $\beta$ which is a linear combination of those.
Then we state and prove \Cref{lem:IncrementPotential} which is our main lemma, and show how it implies \Cref{thm:k=10}.
Lastly, we prove \Cref{lem:PkPm}.
In \Cref{sec:concrems} we discuss further directions and open problems.

\section{Painter's Strategy}
\label{sec:strategy}

Before describing Painter's strategy we introduce some terminology.
We say that a red component is \emph{small} if it is a copy of $P_2$.
Otherwise, we say it is \emph{large}.
When a red component becomes large, three consecutive vertices in it are arbitrarily designated as \emph{central} vertices.
Every vertex of a large red component will be called either central, \emph{outer} or \emph{terminal}.

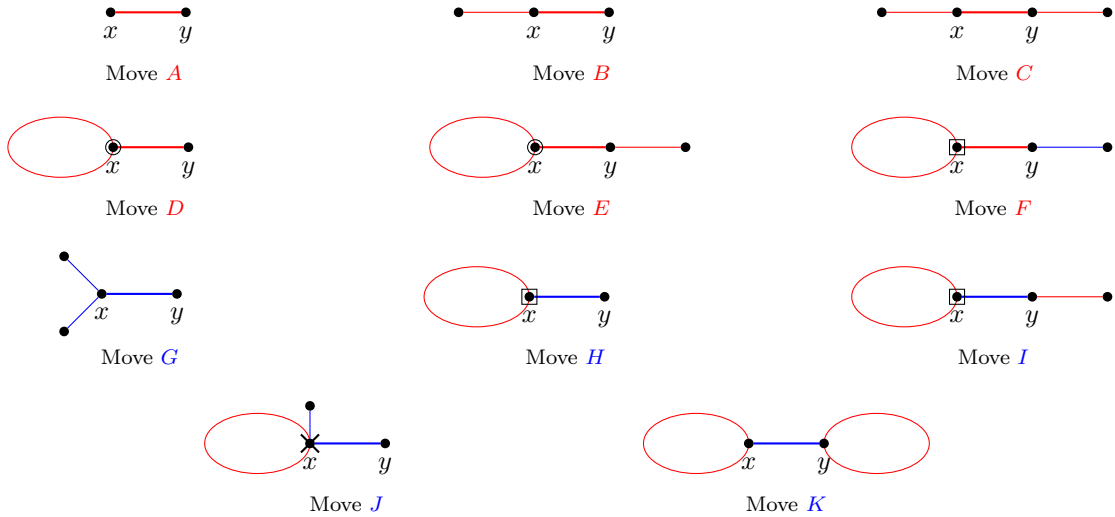
\begin{figure}
\captionsetup[subfigure]{labelformat=empty}
\centering
\subfloat[Move {\color{red}$A$}]{
\begin{tikzpicture}[auto, vertex/.style={circle,draw=black!100,fill=black!100, thick,
			inner sep=0pt,minimum size=1mm},  novertex/.style={circle,draw=white!100}]
		\node (dx) at (-1.7,0) [novertex] {};
        \node (dy) at (1.7,0) [novertex] {};
        \node (x) at (-0.5,0) [vertex,label=below:$x$] {};
		\node (y) at (0.5,0) [vertex,label=below:$y$] {};
		
		\draw [-,red,thick] (x) --node[inner sep=2pt,swap]{} (y);
		\end{tikzpicture}}
\hfill
\subfloat[Move {\color{red}$B$}]{\begin{tikzpicture}[auto, vertex/.style={circle,draw=black!100,fill=black!100, thick,
			inner sep=0pt,minimum size=1mm}, novertex/.style={circle,draw=white!100}]
		\node (dx) at (-1.7,0) [novertex] {};
        \node (dy) at (1.7,0) [novertex] {};
        \node (x) at (-0.5,0) [vertex,label=below:$x$] {};
		\node (y) at (0.5,0) [vertex,label=below:$y$] {};
        \node (x1) at (-1.5,0) [vertex] {};
		
		\draw [-,red,thick] (x) --node[inner sep=2pt,swap]{} (y);
        \draw [-,red] (x1) --node[inner sep=2pt,swap]{} (x);
		\end{tikzpicture}}
\hfill
\subfloat[Move {\color{red}$C$}]{\begin{tikzpicture}[auto, vertex/.style={circle,draw=black!100,fill=black!100, thick,
			inner sep=0pt,minimum size=1mm},  novertex/.style={circle,draw=white!100}]
        \node (dx) at (-1.7,0) [novertex] {};
        \node (dy) at (1.7,0) [novertex] {};
		\node (x) at (-0.5,0) [vertex,label=below:$x$] {};
		\node (y) at (0.5,0) [vertex,label=below:$y$] {};
        \node (x1) at (-1.5,0) [vertex] {};
        \node (y1) at (1.5,0) [vertex] {};
		
		\draw [-,red,thick] (x) --node[inner sep=2pt,swap]{} (y);
        \draw [-,red] (x1) --node[inner sep=2pt,swap]{} (x);
        \draw [-,red] (y1) --node[inner sep=2pt,swap]{} (y);
		\end{tikzpicture}}

\subfloat[Move {\color{red}$D$}]{
\begin{tikzpicture}[auto, vertex/.style={circle,draw=black!100,fill=black!100, thick,
			inner sep=0pt,minimum size=1mm},  novertex/.style={circle,draw=white!100}]
        \node (dy) at (1.7,0) [novertex] {};
		\draw [red] (-1.2,0) ellipse [x radius = 0.7cm, y radius = 0.4cm];
        \node (x) at (-0.5,0) [vertex,label=below:$x$] {};
		\node (y) at (0.5,0) [vertex,label=below:$y$] {};
		
		\draw [-,red,thick] (x) --node[inner sep=2pt,swap]{} (y);

        \draw (-0.5,0) circle [radius = 0.1cm];
        
		\end{tikzpicture}}
\hfill
\subfloat[Move {\color{red}$E$}]{\begin{tikzpicture}[auto, vertex/.style={circle,draw=black!100,fill=black!100, thick,
			inner sep=0pt,minimum size=1mm},  novertex/.style={circle,draw=white!100}]
        \node (dy) at (1.7,0) [novertex] {};
		\draw [red] (-1.2,0) ellipse [x radius = 0.7cm, y radius = 0.4cm];
        \node (x) at (-0.5,0) [vertex,label=below:$x$] {};
		\node (y) at (0.5,0) [vertex,label=below:$y$] {};
        \node (y1) at (1.5,0) [vertex] {};
		
		\draw [-,red,thick] (x) --node[inner sep=2pt,swap]{} (y);
        \draw [-,red] (y1) --node[inner sep=2pt,swap]{} (y);

        \draw (-0.5,0) circle [radius = 0.1cm];
		\end{tikzpicture}}
\hfill
\subfloat[Move {\color{red}$F$}]{\begin{tikzpicture}[auto, vertex/.style={circle,draw=black!100,fill=black!100, thick,
			inner sep=0pt,minimum size=1mm},  novertex/.style={circle,draw=white!100}]
        \node (dy) at (1.7,0) [novertex] {};
		\draw [red] (-1.2,0) ellipse [x radius = 0.7cm, y radius = 0.4cm];
        \node (x) at (-0.5,0) [vertex,label=below:$x$] {};
		\node (y) at (0.5,0) [vertex,label=below:$y$] {};
        \node (y1) at (1.5,0) [vertex] {};
		
		\draw [-,red,thick] (x) --node[inner sep=2pt,swap]{} (y);

        \draw [-,blue] (y1) --node[inner sep=2pt,swap]{} (y);

        \draw (-0.6,-0.1) rectangle +(0.2,0.2);
		\end{tikzpicture}}

\subfloat[Move {\color{blue}$G$}]{\begin{tikzpicture}[auto, vertex/.style={circle,draw=black!100,fill=black!100, thick,
			inner sep=0pt,minimum size=1mm}, novertex/.style={circle,draw=white!100}]
        \node (dx) at (-1.7,0) [novertex] {};
        \node (dy) at (1.7,0) [novertex] {};
        \node (x) at (-0.5,0) [vertex,label=below:$x$] {};
		\node (y) at (0.5,0) [vertex,label=below:$y$] {};
        \node (x1) at (-1,0.5) [vertex] {};
        \node (x2) at (-1,-0.5) [vertex] {};
		
		\draw [-,blue,thick] (x) --node[inner sep=2pt,swap]{} (y);
        \draw [-,blue] (x) --node[inner sep=2pt,swap]{} (x1);
        \draw [-,blue] (x) --node[inner sep=2pt,swap]{} (x2);
		\end{tikzpicture}}
\hfill
\subfloat[Move {\color{blue}$H$}]{\begin{tikzpicture}[auto, vertex/.style={circle,draw=black!100,fill=black!100, thick,
			inner sep=0pt,minimum size=1mm},  novertex/.style={circle,draw=white!100}]
        \node (dy) at (1.7,0) [novertex] {};
		\draw [red] (-1.2,0) ellipse [x radius = 0.7cm, y radius = 0.4cm];
        \node (x) at (-0.5,0) [vertex,label=below:$x$] {};
		\node (y) at (0.5,0) [vertex,label=below:$y$] {};
		
		\draw [-,blue,thick] (x) --node[inner sep=2pt,swap]{} (y);

        \draw (-0.6,-0.1) rectangle +(0.2,0.2);
		\end{tikzpicture}}
\hfill
\subfloat[Move {\color{blue}$I$}]{\begin{tikzpicture}[auto, vertex/.style={circle,draw=black!100,fill=black!100, thick,
			inner sep=0pt,minimum size=1mm},  novertex/.style={circle,draw=white!100}]
        \node (dy) at (1.7,0) [novertex] {};
		\draw [red] (-1.2,0) ellipse [x radius = 0.7cm, y radius = 0.4cm];
        \node (x) at (-0.5,0) [vertex,label=below:$x$] {};
		\node (y) at (0.5,0) [vertex,label=below:$y$] {};
        \node (y1) at (1.5,0) [vertex] {};
		
        \draw [-,red] (y1) --node[inner sep=2pt,swap]{} (y);

        \draw [-,blue,thick] (x) --node[inner sep=2pt,swap]{} (y);

        \draw (-0.6,-0.1) rectangle +(0.2,0.2);
		\end{tikzpicture}}

\subfloat[Move {\color{blue}$J$}]{\begin{tikzpicture}[auto, vertex/.style={circle,draw=black!100,fill=black!100, thick,
			inner sep=0pt,minimum size=1mm}, novertex/.style={circle,draw=white!100}]
        \node (dy) at (1.7,0) [novertex] {};
		\draw [red] (-1.2,0) ellipse [x radius = 0.7cm, y radius = 0.4cm];
        \node (x) at (-0.5,0) [vertex,label=below:$x$] {};
        \draw[thick] (-0.62,0.12) -- (-0.38,-0.12);
        \draw[thick] (-0.62,-0.12) -- (-0.38,0.12);
		\node (y) at (0.5,0) [vertex,label=below:$y$] {};
        \node (x1) at (-0.5,0.5) [vertex] {};
		
		\draw [-,blue,thick] (x) --node[inner sep=2pt,swap]{} (y);
        \draw [-,blue] (x) --node[inner sep=2pt,swap]{} (x1);
		\end{tikzpicture}}
\hspace{55pt}
\subfloat[Move {\color{blue}$K$}]{\begin{tikzpicture}[auto, vertex/.style={circle,draw=black!100,fill=black!100, thick,
			inner sep=0pt,minimum size=1mm}]
		\draw [red] (-1.2,0) ellipse [x radius = 0.7cm, y radius = 0.4cm];
        \draw [red] (1.2,0) ellipse [x radius = 0.7cm, y radius = 0.4cm];
        \node (x) at (-0.5,0) [vertex,label=below:$x$] {};
		\node (y) at (0.5,0) [vertex,label=below:$y$] {};
        
		\draw [-,blue,thick] (x) --node[inner sep=2pt,swap]{} (y);
		\end{tikzpicture}}

\caption{Painter's moves. A central vertex is represented with a circle around it, an outer vertex with a square around it and a terminal vertex is represented with a cross. A simple node represents any vertex.}
\label{fig:moves}
\end{figure}

We present Painter's strategy by spelling out the cases where an edge added by Builder is coloured red and the cases where it is coloured blue.
In principle, it is enough to spell out the cases where an edge is coloured red and let Painter colour it blue otherwise.
However, since we analyse the course of the game based on the set of moves of Painter, and we use them to track the change in values of parameters in $\cP$ (defined in the next section), it will be more convenient to describe both the red and the blue moves explicitly.
Then, we show that this strategy is exhaustive, meaning that it provides Painter with a response for each possible move of Builder.
\begin{strategy}[Strategy for Painter]
\label{strategy:painter}
Let $xy$ be the edge added by Builder in the last round.
In all moves below, $x$ and $y$ are interchangeable whenever their descriptions are not symmetric.
Painter colours $xy$ red in the cases below, always assuming that the previous cases do not hold.
See \Cref{fig:moves} for illustrations of the moves.
\begin{itemize}
    \item[{\color{red}A}] Neither $x$ nor $y$ is adjacent to any red edge.
    \item[{\color{red}B}] $x$ is a vertex in a small red component and $y$ is not adjacent to any red edge.
    Painter calls all three vertices of the new red component central vertices.
    \item[{\color{red}C}] Both $x$ and $y$ are vertices of two different small red components.
    Painter arbitrarily chooses three consecutive vertices on the obtained red $P_4$ and calls them central vertices.
    Painter calls the fourth vertex an outer vertex.
    \item[{\color{red}D}] $x$ is a central vertex of a large red component and $y$ is not adjacent to any red edge.
    Painter calls $y$ an outer vertex.
    \item[{\color{red}E}] $x$ is a central vertex of a large red component and $y$ is a vertex in a small red component.
    Painter calls $y$ and the other vertex of the small red component outer vertices.
    \item[{\color{red}F}] $x$ is an outer vertex of a large component and $y$ is adjacent to precisely one blue edge and no red edges.
    Painter calls $y$ a terminal vertex.
\end{itemize}
Painter colours $xy$ blue in the cases below, always assuming that the previous cases do not hold.
\begin{itemize}
    \item[{\color{blue}G}] Either $x$ or $y$ is incident to at least $2$ blue edges.
    \item[{\color{blue}H}] $x$ is an outer vertex of a large red component and $y$ is not adjacent to any red nor blue edge.
    \item[{\color{blue}I}] $x$ is an outer vertex of a large red component and $y$ is a vertex in a small red component.
    \item[{\color{blue}J}] $x$ is a terminal vertex of a large red component.
    \item[{\color{blue}K}] Both $x$ and $y$ are vertices of large red components.
\end{itemize}
\end{strategy}

\begin{remark}
\label{rem:exhaustive}
    \Cref{strategy:painter} is exhaustive.
    Indeed, if $xy$ is the edge added by Builder in round $t$ of the game, then each of $x$ and $y$ can be exactly one of the following five options: (1) a vertex not contained in any red component; (2) a vertex on a small red component; (3) a central vertex on a large red component; (4) an outer vertex on a large red component; (5) a terminal vertex on a large red component.
    \Cref{table:exhaustive} shows that \Cref{strategy:painter} covers all possible combinations of options for $x$ and $y$.
\end{remark}

\begin{table}[]
\centering
\begin{tabular}{cc|ccccc|}
\cline{3-7}
\multicolumn{2}{c|}{\multirow{2}{*}{}}           & \multicolumn{5}{c|}{$x$}                                                                                        \\ \cline{3-7} 
\multicolumn{2}{c|}{}                            & \multicolumn{1}{c|}{(1)} & \multicolumn{1}{c|}{(2)} & \multicolumn{1}{c|}{(3)} & \multicolumn{1}{c|}{(4)} & (5) \\ \hline
\multicolumn{1}{|c|}{\multirow{5}{*}{$y$}} & (1) & \multicolumn{1}{c|}{{\color{red}$A$}}   & \multicolumn{1}{c|}{{\color{red}$B$}}   & \multicolumn{1}{c|}{{\color{red}$D$}}   & \multicolumn{1}{c|}{{\color{red}$F$},{\color{blue}$G$},{\color{blue}$H$}} & {\color{blue}$G$},{\color{blue}$J$}   \\ \cline{2-7} 
\multicolumn{1}{|c|}{}                     & (2) & \multicolumn{1}{c|}{}    & \multicolumn{1}{c|}{{\color{red}$C$}}   & \multicolumn{1}{c|}{{\color{red}$E$}}   & \multicolumn{1}{c|}{{\color{blue}$G$},{\color{blue}$I$}}   & {\color{blue}$G$},{\color{blue}$J$}   \\ \cline{2-7} 
\multicolumn{1}{|c|}{}                     & (3) & \multicolumn{1}{c|}{}    & \multicolumn{1}{c|}{}    & \multicolumn{1}{c|}{{\color{blue}$G$},{\color{blue}$K$}}   & \multicolumn{1}{c|}{{\color{blue}$G$},{\color{blue}$K$}}   & {\color{blue}$G$},{\color{blue}$J$}   \\ \cline{2-7} 
\multicolumn{1}{|c|}{}                     & (4) & \multicolumn{1}{c|}{}    & \multicolumn{1}{c|}{}    & \multicolumn{1}{c|}{}    & \multicolumn{1}{c|}{{\color{blue}$G$},{\color{blue}$K$}}   & {\color{blue}$G$},{\color{blue}$J$}   \\ \cline{2-7} 
\multicolumn{1}{|c|}{}                     & (5) & \multicolumn{1}{c|}{}    & \multicolumn{1}{c|}{}    & \multicolumn{1}{c|}{}    & \multicolumn{1}{c|}{}    & {\color{blue}$G$},{\color{blue}$J$}   \\ \hline
\end{tabular}
\caption{Painter's moves as responses to any move by Builder, when following \Cref{strategy:painter}. Multiple moves in the same table cell mean that \Cref{strategy:painter} further splits into cases according to the number of incident blue edges to either $x$ or $y$. The meanings of cases (1)-(5) are as in \Cref{rem:exhaustive}.}
\label{table:exhaustive}
\end{table}

We say that a vertex of a red component is of \emph{type 0} if it is not adjacent to any blue edge and of \emph{type 1} if it is adjacent to precisely one blue edge.
\begin{remark}
\label{rem:terminaltype1}
    A terminal vertex is never of type $0$.
    Indeed, a vertex $y$ is labelled terminal only when move $F$ is played, and for this to happen it must be already incident to precisely one blue edge.
\end{remark}

The next lemma captures a relatively simple but important property of \Cref{strategy:painter}.
\begin{lemma}
\label{lem:NoRedPk}
    If Painter follows \Cref{strategy:painter} then no red copy of $P_{10}$ is ever built in the game $\tilde{R}(P_{10},P_n)$.
\end{lemma}

\begin{proof}
We prove this lemma in a series of claims.
For two vertices $x,y \in \mathbb N$ we write $\distR(x,y)$ for the graph distance between them in the red graph.
\begin{Claim}
    Every red component is a tree.
\end{Claim}
\begin{proof}
    By \Cref{strategy:painter}, whenever Painter colours an edge $e=xy$ red we have that $x$ and $y$ are in two different red components.
    Therefore, no cycle inside a red component is ever created.
\end{proof}

\begin{Claim}
    Every large red component contains precisely three central vertices.
    Moreover, if $x$ and $y$ are two central vertices in the same red component, then $\distR(x,y) \leq 2$.
\end{Claim}

\begin{proof}
The only moves which create a large red component are moves $B$ and $C$ and, by definition, for either of them, the new large red component contains precisely three central vertices which are consecutive vertices.
Hence, it suffices to show that the vertices which later join a large red component are never called central vertices.
This is indeed true, as moves $D$, $E$ and $F$ are the only moves where vertices are added to a large red component, and by definition, new added vertices are always called either outer or terminal vertices.
\end{proof}

\begin{Claim}
    If $z$ is an outer vertex in a large red component, then there exists a central vertex $c$ in the same large red component such that $\distR(c,z) \leq 2$.
\end{Claim}
\begin{proof}
    The vertex $z$ was called an outer vertex when Painter played either move $C$, $D$ or $E$.
    For each of those moves it is easy to verify that there exists a central vertex $c$ of the same large red component for which $\distR(c,z) \leq 2$.
\end{proof}

\begin{Claim}
    For every terminal vertex $w$ in a large red component, there exists an outer vertex $z$ in the same large red component for which $\distR(w,z)=1$.
\end{Claim}
\begin{proof}
    This follows easily as $w$ was called terminal when Painter played move $F$, and therefore shares a red edge with an outer vertex $z$.
\end{proof}

We are now ready to finish the proof. Let $w_1$ and $w_2$ be two vertices of a red component. Suppose first that $w_1$ and $w_2$ are both terminal vertices.
Then according to the claims above, there exist outer vertices $z_1$ and $z_2$ and central vertices $c_1$ and $c_2$ of the same large red component, such that $\distR(w_1,z_1)=1$, $\distR(w_2,z_2)=1$, $\distR(z_1,c_1)\leq 2$ and $\distR(z_2,c_2)\leq 2$.
Moreover, we know that $\distR(c_1,c_2) \leq 2$.
Therefore by the triangle inequality we have 
\begin{align*}
    \distR(w_1,w_2) &\leq \distR(w_1,z_1) + \distR(z_1,c_1) + \distR(c_1,c_2)+\distR(c_2,z_2) + \distR(z_2,w_2) \\
    &\leq 1+2+2+2+1 \\
    &=8.
\end{align*}
The other cases can be analysed the same way, and allow us to conclude that $\distR(w_1,w_2) \leq 8$ for any two vertices $w_1$ and $w_2$ in the same large red component.
Thus each red component is a tree of diameter at most $8$, and therefore does not contain any copy of $P_{10}$.
\end{proof}

\section{Potential Function and Main Lemma}
\label{sec:potentialfunc}

As explained in \Cref{sec:idea}, the potential function $\beta$ is derived from the solution of the linear programming problem associated to \Cref{strategy:painter}.
We analyse the course of the game, assuming that Painter follows \Cref{strategy:painter}, by tracking the values of the following parameters of the board of the game.
\begin{itemize}
    \item $\Dn$: Number of vertices of degree one in the blue graph.
    \item $\Dt$: Number of vertices of degree at least two in the blue graph.
    \item $\So$: Number of vertices of type $0$ on small red components.
    \item $\Sn$: Number of vertices of type $1$ on small red components.
    \item $\Co$: Number of central vertices of type $0$.
    \item $\Cn$: Number of central vertices of type $1$.
    \item $\Oo$: Number of outer vertices of type $0$.
    \item $\On$: Number of outer vertices of type $1$.
    \item $\Tr$: Number of terminal vertices of type $1$.
\end{itemize}
As mentioned in \Cref{sec:idea}, we denote this set of parameters by $\cP$.
Moreover, this set of parameters is in fact a set of functions from $\mathbb Z_{\ge 0}$ to $\mathbb Z_{\ge 0}$ (in particular, non-negative), where for an integer $t \in \mathbb N$ and for $X \in \cP$, we write $X(t)$ for the value of the parameter $X$ for the board of the game after precisely $t$ rounds have been completed.

We are now ready to define our potential function.
\begin{equation}
\label{eq:beta}
    \beta(t) \coloneqq  \frac{2}{3}\Dn(t) + \frac{5}{3}\Dt(t) + \frac{1}{2}\So(t) + \frac{1}{2}\Sn(t) + \frac{2}{3}\Co(t) + \frac{1}{2}\Cn(t) + \Oo(t) + \frac{2}{3}\On(t) + \Tr(t).
\end{equation}

\Cref{thm:k=10} follows easily from the next lemma, which is the key ingredient on our proof of it.
\begin{lemma}
\label{lem:IncrementPotential}
    We have $\beta(0) = 0$.
    Moreover, if Painter follows \Cref{strategy:painter} and the game $\tilde{R}(P_{10},P_n)$ lasts for $N$ rounds, for some integer $N \in \mathbb N$, then for every $t \in [N]$ we have
    \[\beta(t) - \beta(t-1) \le 1. \]
\end{lemma}

Before proving \Cref{lem:IncrementPotential}, we show how it readily implies \Cref{thm:k=10}.
\begin{proof}[Proof of \Cref{thm:k=10} assuming \Cref{lem:IncrementPotential}.]
    Suppose that Painter follows \Cref{strategy:painter} and the game $\tilde{R}(P_{10},P_n)$ lasts for $N$ rounds, for some integer $N \in \mathbb N$.
    By \Cref{lem:NoRedPk}, when the game ends there is a blue copy of $P_n$, so we clearly have $\Dn(N) + \Dt(N) \ge n$ and $\Dt(N) \ge n-2$.
    Recall that all parameters in $\cP$ have non-negative values, so we get
    \[\beta(N) \ge \frac{2}{3}\Dn(N) + \frac{5}{3}\Dt(N) \ge \frac{2}{3}\left(n - \Dt(N) \right) + \frac{5}{3}\Dt(N) \ge \frac{5}{3}n - 2. \]
    By \Cref{lem:IncrementPotential} we get that
    \[N \ge \sum_{t\in [N]}\beta(t) - \beta(t-1) = \beta(N). \]
    Combining both inequalities above we get $N \ge \frac{5}{3}n-2$, as we wanted.
\end{proof}

It is left to prove \Cref{lem:IncrementPotential}.
For any function $f : \mathbb{N} \rightarrow \mathbb R$ and $t \in \mathbb N$ we write $\Delta f(t) \coloneqq f(t)-f(t-1)$.

\begin{proof}[Proof of \Cref{lem:IncrementPotential}]
Note that all parameters $X \in \cP$ trivially satisfy $X(0) = 0$, so in particular we have $\beta(0)=0$.
Hence, it is left to show that when Painter follows \Cref{strategy:painter} we have $\Delta \beta(t) \leq 1$ for every round $t \in [N]$ in the game.
To do so, we go over all moves $A$ to $K$, and for each one we analyse the change it causes in parameters in $\cP$ and consequentially in $\beta$, if played by Painter on round $t$ of the game.
We show that all cases lead to $\Delta\beta(t) \leq 1$.
For some of those moves, we will need to distinguish between several cases, depending on the number of blue incident edges to either $x$ or $y$.
However, the case analysis remains elementary.

\begin{itemize}
    \item[{\color{red}$A$}] Let $\alpha_0$ and $\alpha_1$ be the number of vertices of type $0$ and $1$, respectively, amongst $x$ and $y$ at the end of round $t-1$.
    Then we have $\alpha_0 + \alpha_1 \le 2$, $\Delta\So(t)=\alpha_0$, $\Delta\Sn(t)=\alpha_1$, and the other parameters do not change.
    Therefore
    \begin{align*}
    \Delta\beta(t) &= \frac{1}{2}\alpha_0 + \frac{1}{2}\alpha_1 \\
    &= 1.
    \end{align*}
    
    \item[{\color{red}$B$}] Let $x'$ be the vertex such that the edge $x'x$ was red at the end of round $t-1$. 
    Suppose further that at the end of round $t-1$ there were $\alpha_0$ and $\alpha_1$ vertices of type $0$ and $1$, respectively, amongst $x,x'$, and that the number of blue edges incident to $y$ was precisely $a$.
    We have $\alpha_0 \leq \alpha_0 + \alpha_1 \leq 2$, $\Delta\So(t)=-\alpha_0$, $\Delta\Sn(t)=-\alpha_1$, $\Delta\Co(t)=\alpha_0+1_{a=0}$, $\Delta\Cn(t)=\alpha_1+1_{a=1}$, and the other parameters do not change.
    Therefore
    \begin{align*}
    \Delta\beta(t) &=-\frac{1}{2}\alpha_0 + \frac{2}{3}(\alpha_0+1_{a=0}) - \frac{1}{2}\alpha_1 + \frac{1}{2}(\alpha_1+1_{a=1}) \\
    &\leq \frac{1}{6}\alpha_0 + \frac{2}{3}(1_{a=0}+1_{a=1}) \\
    &\leq \frac{1}{6}\cdot 2+\frac{2}{3} \\
    &= 1.
    \end{align*}
    
    \item[{\color{red}$C$}] Let $x',y'$ be the vertices such that the edges $x'x$ and $y'y$ were red at the end of round $t-1$.
    Without loss of generality we may suppose that Painter calls $x$, $x'$ and $y$ central vertices.
    Let $\alpha_0$ and $\alpha_1$ be the number of vertices of type $0$ and $1$, respectively, amongst $x',x,y$ at the end of round $t-1$, so we have $\alpha_0 \le \alpha_0 + \alpha_1 \le 3$.
    Let $a$ be the number of blue edges incident to $y'$ at the end of round $t-1$.
    Then we have $\Delta\So(t)=-\alpha_0-1_{a=0}$, $\Delta\Sn(t)=-\alpha_1-1_{a=1}$, $\Delta\Co(t)=\alpha_0$, $\Delta\Cn(t)=\alpha_1$, $\Delta\Oo(t) = 1_{a=0}$, $\Delta\On(t) = 1_{a=1}$, and the other parameters do not change.
    Therefore
    \begin{align*}
        \Delta\beta(t) &= \frac{1}{2}(-\alpha_0-1_{a=0}) + \frac{1}{2}(-\alpha_1-1_{a=1}) + \frac{2}{3}\alpha_0 + \frac{1}{2}\alpha_1 + 1_{a=0} + \frac{2}{3}\cdot 1_{a=1} \\
        &\leq \frac{1}{6}\alpha_0 +\frac{1}{2}(1_{a=0}+1_{a=1}) \\
        &\leq \frac{1}{6}\cdot 3 +\frac{1}{2} \\
        &= 1.
    \end{align*}
    
    \item[{\color{red}$D$}] Let $a$ be the number of blue edges incident to $y$ at the end of round $t-1$.
    Then we have $\Delta\Oo(t) = 1_{a=0}$, $\Delta\On(t) = 1_{a=1}$, and the other parameters do not change.
    Therefore
    \begin{align*}
    \Delta\beta(t) &= 1_{a=0}+\frac{2}{3}\cdot 1_{a=1} \\
    &\leq 1_{a=0}+1_{a=1} \\
    &\leq 1.
    \end{align*}
    
    \item[{\color{red}$E$}] Let $y'$ be the vertex such that the edge $y'y$ was red at the end of round $t-1$. 
    Let $\alpha_0$ and $\alpha_1$ be the number of vertices of type $0$ and $1$, respectively, amongst $y$ and $y'$.
    We have $\alpha_0 + \alpha_1 \le 2$, $\Delta\So(t)=-\alpha_0$, $\Delta\Sn(t)=-\alpha_1$, $\Delta\Oo(t)=\alpha_0$, $\Delta\On(t)=\alpha_1$, and the other parameters do not change.
    Therefore
    \begin{align*}
    \Delta\beta(t) &= -\frac{1}{2}\alpha_0 + - \frac{1}{2}\alpha_1 + \alpha_0 + \frac{2}{3}\alpha_1 \\
    &\leq \frac{1}{2}(\alpha_0+\alpha_1) \\
    &\leq 1.
    \end{align*}
    
    \item[{\color{red}$F$}] We have $\Delta\Tr(t)=1$, and the other parameters do not change.
    Therefore $\Delta\beta(t) = 1$.

    \item[{\color{blue}$G$}] Define the following function
    \begin{align}
    \label{eq:gamma}
        \gamma(t) \coloneqq \frac{1}{2}\So(t) + \frac{1}{2}\Sn(t) + \frac{2}{3}\Co(t) + \frac{1}{2}\Cn(t) + \Oo(t) + \frac{2}{3}\On(t) + \Tr(t).
    \end{align}
    Then we have $\beta(t) = \frac{2}{3}\Dn(t) + \frac{5}{3}\Dt(t) + \gamma(t)$.

    Assume without loss of generality that, at the end of round $t-1$, the vertex $x$ is adjacent to at least two blue edges and that the number of blue edges incident to $y$ is precisely $a$.
    There are several cases to consider.
    \begin{itemize}
        \item If $y$ is a vertex on a small red component, then we have $\Delta\So(t) = -1_{a=0}$, $\Delta\Sn(t) = 1_{a=0} - 1_{a=1}$ and therefore $\Delta\gamma(t) = -\frac{1}{2}1_{a=0} + \frac{1}{2}(1_{a=0} - 1_{a=1}) \le 0$.
        \item If $y$ is a central vertex, then we have $\Delta\Co(t) = -1_{a=0}$, $\Delta\Cn(t) = 1_{a=0} - 1_{a=1}$ and therefore $\Delta\gamma(t) = -\frac{2}{3}1_{a=0} + \frac{1}{2}(1_{a=0} - 1_{a=1}) \le 0$.
        \item If $y$ is an outer vertex, then we have $\Delta\Oo(t) = -1_{a=0}$, $\Delta\On(t) = 1_{a=0} - 1_{a=1}$ and therefore $\Delta\gamma(t) = -1_{a=0} + \frac{2}{3}(1_{a=0} - 1_{a=1}) \le 0$.
        \item If $y$ is a terminal vertex then we have $\Delta\Tr(t) \le 0$ and therefore $\Delta\gamma(t) = 0$.
        \item If $y$ is neither of the above then clearly we have $\Delta\gamma(t) = 0$.
    \end{itemize}
    In each of these cases we get that $\Delta\gamma(t) \le 0$.
    Furthermore, we find that $\Delta\Dn(t) = 1_{a=0} - 1_{a=1}$ and $\Delta\Dt(t) = 1_{a=1}$.
    Therefore $\Delta\beta(t) \le \frac{2}{3}\cdot 1_{a=0} + 1_{a=1} \le 1$.

    \item[{\color{blue}$H$}] Let $a$ be the number of blue edges incident to $x$ at the end of round $t-1$.
    We then have $\Delta\Dn(t) = 1 + 1_{a=0} - 1_{a=1}$, $\Delta\Dt(t) = 1_{a=1}$, $\Delta\Oo(t) = -1_{a=0}$ and $\Delta\On(t) = 1_{a=0}-1_{a=1}$, and the other parameters do not change.
    Therefore
    \begin{align*}
    \Delta\beta(t) &\le \frac{2}{3}(1+1_{a=0}-1_{a=1}) + \frac{5}{3}\cdot 1_{a=1} - 1_{a=0} + \frac{2}{3}(1_{a=0}-1_{a=1}) \\
    &= \frac{2}{3} + \frac{1}{3}\cdot (1_{a=0} + 1_{a=1}) \\
    &= 1.
    \end{align*}

    \item[{\color{blue}$I$}] Let $a$ and $b$ be the numbers of blue edges incident to $x$ and $y$, respectively, at the end of round $t-1$.
    We then have $\Delta\Dn(t) = 1_{a=0}-1_{a=1}+1_{b=0}-1_{b=1}$, $\Delta\Dt(t) = 1_{a=1}+1_{b=1}$, $\Delta\Oo(t) = -1_{a=0}$, $\Delta\On(t) = 1_{a=0}-1_{a=1}$, $\Delta\So(t)=-1_{b=0}$ and $\Delta\Sn(t) = 1_{b=0}-1_{b=1}$, and the other parameters do not change.
    Therefore 
    \begin{align*}
    \Delta\beta(t) \leq& \frac{2}{3}(1_{a=0}-1_{a=1}+1_{b=0}-1_{b=1}) + \frac{5}{3}(1_{a=1}+1_{b=1}) \\
    &-1_{a=0} +\frac{2}{3}(1_{a=0}-1_{a=1})-\frac{1}{2}\cdot 1_{b=0}+\frac{1}{2}(1_{b=0}-1_{b=1}) \\
    =& \frac{1}{3}(1_{a=0} + 1_{a=1}) + \frac{2}{3}\cdot 1_{b=0} + \frac{1}{2}\cdot 1_{b=1} \\
    \leq& 1.
    \end{align*}

    \item[{\color{blue}$J$}] Since move $G$ was not played, the terminal vertex $x$ is adjacent to at most one blue edge, and by \Cref{rem:terminaltype1}, $x$ is adjacent to precisely one blue edge.
    In particular, we get $\Delta\Tr(t) \le -1$.
    Considering this, and repeating the case analysis from move $G$ for the value of the function $\gamma$ as defined in (\ref{eq:gamma}), we obtain $\Delta\gamma(t) \le -1$.
    Let $a$ be the number of blue edges incident to $y$ at the end of round $t-1$.
    Hence, we also have $\Delta\Dn(t)= -1 -1_{a=1}+1_{a=0} = -2\cdot 1_{a=1}$, $\Delta\Dt(t)= 1+1_{a=1}$, and therefore
    \begin{align*}
    \Delta\beta(t) &\le -\frac{2}{3}\cdot 2 \cdot 1_{a=1} + \frac{5}{3}(1+1_{a=1}) - 1 \\
    &\le \frac{2}{3} + \frac{1}{3}\cdot 1_{a=1} \\
    &\le 1.
    \end{align*}

    \item[{\color{blue}$K$}] Let $c_0$, $c_1$, $o_0$, $o_1$ be respectively the numbers of central vertices of type $0$, of central vertices of type $1$, of outer vertices of type $0$ and of outer vertices of type $1$ amongst $x,y$ at the end of round $t-1$.
    Then we have $c_0+c_1+o_0+o_1=2$, $\Delta\Dn(t) = c_0+o_0-c_1-o_1$, $\Delta\Dt(t) = c_1+o_1$, $\Delta\Co(t)=-c_0$, $\Delta\Cn(t)=c_0-c_1$, $\Delta\Oo(t)=-o_0$, $\Delta\On(t)=o_0-o_1$, and the other parameters do not change.
    Therefore
    \begin{align*}
    \Delta\beta(t) &= \frac{2}{3}(c_0+o_0-c_1-o_1) + \frac{5}{3}(c_1+o_1) -\frac{2}{3}c_0 + \frac{1}{2}(c_0-c_1) - o_0 + \frac{2}{3}(o_0-o_1) \\
    &= \frac{1}{2}c_0 + \frac{1}{2}c_1 + \frac{1}{3}o_0 + \frac{1}{3}o_1 \\
    &\le \frac{1}{2}(c_0+c_1+o_0+o_1)\\
    &= 1.
    \end{align*}
\end{itemize}
Therefore, we get that $\Delta\beta(t) \le 1$ for every $t\in [N]$, proving \Cref{lem:IncrementPotential} and so \Cref{thm:k=10}.
\end{proof}

It is then left to prove \Cref{lem:PkPm}, which is rather simple.
\begin{proof}
    Consider the game $\tilde{R}(P_k,H)$.
    Painter follows a strategy that guarantees that neither a red copy of $P_m$ nor a blue copy of $H$ is created for $\tilde{r}(P_m,H)-1$ many rounds.
    In particular, all red paths at this point of the game are of length at most $m-1$.
    Then, from the $\tilde{r}(P_m,H)$-th round and on, Painter colours red any edge added by Builder.
    It is easy to see that then Builder wins once a red copy of $P_k$ is created.
    Suppose it takes Builder $\tilde{r}(P_m,H)-1+\ell$ many rounds to win, for some integer $\ell \ge 1$, and let $P$ be a red copy of $P_k$ at the end of the game.
    Partition the $k-1$ edges of $P$ into $\left\lceil\frac{k-1}{m-1} \right\rceil$ many parts, each but at most one consists of $m-1$ consecutive edges.
    Each of these $\left\lfloor\frac{k-1}{m-1} \right\rfloor$ parts is, in fact, a red copy of $P_m$.
    Hence, in each one of them, at least one edge was added by Builder during the last $\ell$ rounds of the game, implying $\ell \ge \left\lfloor\frac{k-1}{m-1} \right\rfloor$ and finishing the proof.
\end{proof}

\section{Concluding remarks}
\label{sec:concrems}

As we mentioned in \Cref{cor:main}, we determine the asymptotic value of $\tilde{r}(P_k,P_n)$ for every $10 \le k = o(n)$, by matching our lower bound from \Cref{thm:main} with the upper bound from~\cite{BEDNARSKABZDEGA2024103873}.
In the case $k=3$, Cyman, Dzido, Lapinskas and Lo~\cite{cyman2015line} showed that $\tilde{r}(P_3,P_n)= \left\lceil \frac{5(n-1)}{4} \right\rceil$.
For $k=4$, they also showed that $\tilde{r}(P_4,P_n) \geq \left\lceil \frac{7n}{5}-1 \right\rceil$, and later it was proved independently by Bednarska-Bzd{\k{e}}ga~\cite{BEDNARSKABZDEGA2024103873} and Zhang and Zhang~\cite{zhang2023proof} that in fact $\tilde{r}(P_4,P_n) = \left\lceil \frac{7n}{5}-1 \right\rceil$.

The most natural continuation of this work would be to obtain exact values of $\tilde{r}(P_k,P_n)$ for all values of $k$ and $n$.
However, for $5 \le k \le 9$ or when $k = \Theta(n)$, even the asymptotic value of $\tilde{r}(P_k,P_n)$ is still unknown, or whether the limit exists.

While in this paper we considered paths, another natural line of research would be to study online Ramsey numbers for cycles (see \cite{adamski2023online, adamski2021online}), or of cycles and paths (see \cite{adamski2022online}), which may be closely related to it, as shown in~\cite{cyman2015line}.

\paragraph*{Acknowledgements.}
The authors would like to thank their PhD supervisor Professor B\'{e}la Bollob\'{a}s for his support and valuable comments.
The authors are grateful to Natalia Adamska for pointing out an error in a previous version of this manuscript. The authors would also like to thank Grzegorz Adamski and Ma\l gorzata Bednarska-Bzd{\k{e}}ga for pointing out that the proof of $\tilde{r}(P_k,P_n) \geq \frac{5}{3} n -2$ in a previous version of this manuscript actually yields the stronger bound in \Cref{thm:main} by a simple application of their suggested \Cref{lem:PkPm}. \\
The first author is funded by Trinity College, Cambridge.
The second author is funded by EPSRC (Engineering and Physical Sciences Research Council) and by the Cambridge Commonwealth, European and International Trust.

\bibliographystyle{amsplain}  
\renewcommand{\bibname}{Bibliography}
\bibliography{bib}

\end{document}